\documentclass[a4paper,12pt]{amsart}
\usepackage{amssymb}
\usepackage{amsopn} 
\usepackage{amsthm}
\usepackage{amsmath}
\newtheorem{lemma}{Lemma}[section]
\newtheorem{prop}{Proposition}[section]
\newtheorem{definition}{Definition}[section]
\newtheorem{introtheorem}{Theorem}

\newtheorem{theorem}{Theorem}
\newtheorem{corollary}[theorem]{Corollary}

\newcommand{\GmodG}{{\G_1\!\backslash\G}}
\newcommand{\N}{{\mathbb N}}

\newcommand{\R}{{\mathbb R}}
\newcommand{\C}{{\mathbb C}}

\newcommand{\w}{{\omega}}
\newcommand{\SL}{\operatorname{SL_2}}

\newcommand{\abs}[1]{\left\lvert #1 \right\rvert}
\newcommand{\norm}[1]{\left\lVert #1 \right\rVert}
\newcommand{\mymatrix}[4]{\left(\begin{array}{cc}#1&#2\\#3&#4\end{array}\right)}
\newcommand{\modsym}[2]{\left\langle #1,#2 \right\rangle}
\newcommand{\inprod}[2]{\left \langle #1,#2 \right\rangle}

\DeclareMathOperator*{\res}{Res}
\DeclareMathOperator*{\ord}{ord}
\DeclareMathOperator*{\vol}{vol}

\def \GmodH {{\Gamma\!\setminus\!\H}}
\def \G {{\Gamma}}
\def \g {{\gamma}}
\def \a {{\alpha}}
\def \d {{\delta}}
\def \la {{\lambda}}
\def \L {{\operatorname{L}}}
\def \s {{\sigma}}

\def \e {{\epsilon}}
\def \ev {{\vec\epsilon}}
\def \ov {{\vec{0}}}
\def \H {{\mathbb H}}

\def \psl {{\operatorname{PSL}_2(\mathbb{R})}}

\title{Distribution of modular symbols for compact surfaces} 

\begin{document}
\author{Morten Skarsholm Risager$^\dagger$}
\address{Department of Mathematical Sciences\\ University of Aarhus\\ Ny Munkegade Building 530\\ 8000 {Aa}rhus, Denmark}
\email{risager@imf.au.dk}
\date{July 31, 2003}
\thanks{$^\dagger$Partially supported by MaPhySto - A Network in Mathematical Physics and Stochastics, funded by The Danish National Research Foundation.
}
\subjclass[2000]{11F67, 11F72, 11M36}
\begin{abstract}
We prove that the modular symbols appropriately normalized and ordered have an asymptotical normal distribution
for all cocompact subgroups of $\hbox{SL}_2( {\mathbb R})$. We introduce hyperbolic Eisenstein series in order to calculate the moments of the modular symbols.  
\end{abstract}
\maketitle
\section{Introduction}
Let $M$ be a hyperbolic Riemann surface of finite volume. % Se Forster 27.10 s.211 for hyperbolic
  Hence the universal covering of $M$ is the upper half plane, $\H$, and the covering group, $\G$, is a discrete subgroup of $\psl$.
 Let $f(z)dz$ be a holomorphic 1-form on $M$. If $c$ is a curve on $M$ we may  integrate $f(z)dz$ along the curve to get 
$$\int_cf(z)dz.$$

We have a bijection between the covering group $\G$ and the fundamental group $\pi_1(M,\hat z_0)$ %See forster Theorem 5.6
given by sending $\gamma\in \G$ to the unique geodesic between $z_0$ and $\gamma z_0$ in $\H$ where $z_0$ lies above $\hat z_0$, and then projecting this curve to $M$. By integrating along this curve we get an additive homomorphism
\begin{align}\label{map}\G&\to\qquad\C\\
                      \nonumber\g&\mapsto\,\int_{z_0}^{\g z_0}\!\!\!f(z)dz.\end{align}
We wish to study the distribution of the values of this  map.

In \cite{petridisrisager} we considered this map in the case where $M$ is non-compact and  $f$ is cuspidal. In applications to analytic number theory and elliptic curves this is often the relevant setup. In topology and geometry, on the other hand, the case of $M$ compact usually attracts more attention.  In  \cite[Theorem B]{petridisrisager} we found that with the correct normalization and ordering the values of the map (\ref{map}) are normally distributed when $M$ has a cusp and $f$ is cuspidal. In this paper we obtain a similar result in the case where $M$ is compact and  $f$ satisfies a different condition  than being cuspidal.

 The result in  \cite{petridisrisager} was proved using Eisenstein series twisted with modular symbols introduced by Goldfeld (\cite{goldfeld2,goldfeld1}). The definition of these requires the existence of a cusp and can therefore \emph{not} be used in the compact case. In this paper we introduce hyperbolic Eisenstein series.  They exist also in the compact case. We then \lq twist\rq{} these  with modular symbols to obtain the distributional result described below. 

We can handle a slightly more general setup than described above. Let $\G\subseteq\psl$ be discrete and cocompact and  let $M=\GmodH$ be the associated quotient space. Let $\g_1\in \G$ be hyperbolic, i.e. $\abs{\operatorname{tr}(\g_1)}>2$. For simplicity we assume $$\g_1=\mymatrix{\sqrt{\mu}}{0}{0}{\sqrt{\mu}^{-1}} $$ where $\mu>0$. This may always be obtained from any hyperbolic $\g_1$ by conjugation with $g\in\SL(\R)$. By possibly considering $\g_1^{-1}$ instead of $\g_1$ we may assume $\mu>1$. We further assume that $f(z)dz$ is a holomorphic 1-form on $\GmodH$  which satisfies  $$\int_{z_0}^{\g_1z_0}f(z)dz=0.$$ (We note that such $f$ always exist whenever the genus of $\GmodH$ is strictly larger than 1). We define 
 \begin{align*}[\g,f]=\sqrt{\frac{\vol{(\GmodH)}}{\log((a^2+b^2)(c^2+d^2))}}\int_{z_0}^{\g z_0}\!\!\!f(z)dz,\end{align*} which for fixed $f$ gives a map from the quotient  $\GmodG$ to $\C$.
Here $\G_1$ is the cyclic subgroup of $\G$ generated by $\g_1$ and $a,b,c$ and $d$ are the entries of $\g$.

Our main theorem is the following distributional result.
\begin{introtheorem}\label{maintheoremforcomplex}    Assume $f$ has Petersson norm 1. Then $[\g,f]$ has asymptotical normal distribution. 
More precisely, for any fixed rectangle $R$ in $\C$, 
\begin{equation}\label{mainlimit}\frac{ \#\left\{\gamma\in (\GmodG)^T \left| [\g,f]\in R\right. \right\} }{\#(\GmodG)^T}\to \frac{1}{2\pi}\int_R\!\! \exp\left(-\frac{x^2+y^2}{2}\right)dxdy\end{equation} as $T\to\infty$.
\end{introtheorem}
Here 
$$ (\GmodG)^T =\left\{\g\in\GmodG|\quad \sqrt{(a^2+b^2)(c^2+d^2)}\leq T\right\}.$$

We also have a distribution result for real harmonic differentials $\a=\Re(f(z)dz)$. In this case we find the following:
\begin{introtheorem}\label{maintheoremforreal} Assume $f$ has Petersson norm 1. Then $$[\g,\a]=\sqrt{\frac{\vol{(\GmodH)}}{\log((a^2+b^2)(c^2+d^2))}}\int_{z_0}^{\g z_0}\!\!\!\a$$ has asymptotical normal distribution. 
More precisely, for any fixed rectangle $R$ in $\C$, 
\begin{equation}\label{secondmainlimit}\frac{ \#\left\{\gamma\in (\GmodG)^T \left| [\g,\a]\in R\right. \right\} }{\#(\GmodG)^T}\to \frac{1}{\sqrt{2\pi}}\int_R\!\! \exp\left(-\frac{x^2}{2}\right)dxdy\end{equation} as $T\to\infty$.
\end{introtheorem}

In order to prove such  results we introduce hyperbolic Eisenstein series defined by 
$$ E^{\g_1}(z,s)=\sum_{\g\in\GmodG}\left( \frac{\Im(\g z)}{\abs{\g z}}\right)^s\qquad\textrm{ for } \Re(s)>1.$$
This converges absolutely for  $\Re(s)>1$ by Lemma \ref{convergence} below. We then go on to study the basic properties of this series.
\begin{introtheorem}\label{imgoingtoaconcert} The function $E^{\g_1}(z,s)$ has meromorphic continuation to the whole $s$ plane. At a regular point, $s_0$,  $E^{\g_1}(z,s_0)$ is square integrable on $\GmodH$.  The poles are located at $-2n+s_j$ where $s_j(1-s_j)$ is an eigenvalue of the automorphic Laplacian and $n\in\N$. The point $s=1$ is a simple pole and the  residue at $s=1$ is $$\frac{2\log \mu}{\vol{(\GmodH)}}.$$ For fixed $\s=\Re(s)$, $1/2<\s\leq 1$, the hyperbolic Eisenstein series has at most polynomial growth on vertical lines $\s+it.$
\end{introtheorem}
Most of this follows quite straightforward after applying the resolvent to the following identity
\begin{equation}(\Delta+s(1-s)) E^{\g_1}(z,s)=-s^2 E^{\g_1}(z,s+2).\end{equation}
Once the above theorem is established we can use the method of complex contour integration to get the following result which may be interpreted as a result on the number of closed homotopy classes on a compact hyperbolic Riemann surface.
\begin{introtheorem}\label{homotopyclasses}
\begin{equation*}\sum_{\substack{\gamma\in\GmodG\\\sqrt{(a^2+b^2)(c^2+d^2)}\leq T}}1=\frac{2\log\mu}{\vol(\GmodH)}T+O(T^{1-\delta})\end{equation*} for some $\delta>0$.  \end{introtheorem}
We show that $1-\delta=7/8+\varepsilon$ is valid if there are no small eigenvalues.

Let $f_i$ be modular forms of weight 2 with respect to $\G$ and let $\a_i=\Re(f_i(z)dz)$ or $\a_i=\Im(f_i(z)dz)$. We shall write $\a=\a_1$. The (real) modular symbols are defined by
$$\modsym{\g}{\a_i}=-2\pi i\int_{z_0}^{\g z_0}\a_i.$$ Assume that $\modsym{\g_1}{\a_i}=0$ for $i=1\ldots n$.
  We now ``twist'' the hyperbolic Eisenstein series with these modular symbols as done by Goldfeld \cite{goldfeld2,goldfeld1} for the usual non-holomorphic Eisenstein series by setting
$$E^{\g_1,\a_1,\ldots,\a_n}(z,s)=\sum_{\gamma\in\GmodG} \prod_{k=1}^n\modsym{\g}{\a_k} \left( \frac{\Im(\g z)}{\abs{\g z}}\right)^s\qquad \textrm{ for }\Re(s)\gg 0$$

We then go on to study the analytic properties of this function. 
\begin{introtheorem}\label{whattodo} The function $E^{\g_1,\a_1,\ldots,\a_n}(z,s)$ has meromorphic continuation to the whole $s$-plane. In $\Re(s)>1$ it is  analytic.
\end{introtheorem}
The last claim of the theorem enables us to give rather good bounds on the growth of the modular symbols.
\begin{introtheorem}\label{rathergoodbound} For any $\varepsilon>0$ we have 
$$\modsym{\g}{\a}=O_\varepsilon(((a^2+b^2)(c^2+d^2))^\varepsilon)$$
\end{introtheorem}
We go on to study the possible singularity at  $s=1$. We estimate the pole order and determine the leading term in the Laurent expansion for many cases. We then go on to study the behavior on vertical lines and we arrive at the following thorem.
\begin{introtheorem}\label{lasttheorem}The function $E^{\g_1,\a_1,\ldots,\a_n}(z,s)$ grows at most polynomially on vertical lines with $\Re(s)>1/2$.
\end{introtheorem}  
This puts us in a position where we can use the method of contour integration to calculate the moments of the random variable defined by the left hand side
 of (\ref{mainlimit}). Once we have calculated these moments Theorem \ref{maintheoremforcomplex} follows from a classical theorem in probability theory.  

\hbox{}

{\bf Acknowledgments:}
{I am grateful to Professor A. B. Venkov for drawing my attention to \cite{kudlamillson} and for stimulating discussions regarding hyperbolic Eisenstein series. I am also grateful to Professors E. Balslev and Y. N. Petridis for remarks concerning an early draft of this paper.}

\section{The resolvent of the automorphic Laplacian}\label{resolvent}
For the methods used in this paper it is very important to introduce the resolvent of the automorphic Laplacian. The automorphic Laplacian is closely related to the ordinary hyperbolic Laplacian $$\Delta=y^2\left(\frac{\partial^2}{\partial x^ 2}+\frac{\partial^2}{\partial y^ 2}\right).$$ 
We shall briefly recall the relevant definitions and properties.

Let $\G\subseteq\psl$ be discrete and cocompact and $M=\GmodH$ the associated quotient space under the action
$$\begin{array}{rccc}
\g:&\H&\to&\H\\
   &z&\mapsto&\frac{az+b}{cz+d}.
\end{array}$$
The quotient can be given a structure of a Riemann surface with $\H$ as a branched cover. (See e.g. \cite[\S 1.5]{MR47:3318}) The branch points are at the elliptic points i.e. the points which are fixpoints of $\g\in\G$ with $\abs{\operatorname{tr}(\g)}<2$. When there are no such points $\H$ is the universal cover.

The automorphic Laplacian, $L_\G$ is the closure of the operator acting on smooth forms in $\L^2(\GmodH)$ by $\Delta f$ where $f:\H\to\H$ is $\G$-automorphic and smooth. The operator $L_\G$ is selfadjoint with $-L_\G$ non-negative. By the  maximum principle $L_\G u=0$ if and only if $u$ is constant.  There is a complete orthonormal system of smooth eigenfunctions $\psi_0,\ldots,\psi_i,\ldots$ in $\L^2(\GmodH)$ with 
\begin{align*}
-L_\G\psi_i&=\la_i\psi_i\\
0=\la_0<&\la_1\leq\la_2\leq\ldots
\end{align*}
and $\la_n\to\infty$. (See e.g. \cite[Theorem 7.2.6]{buser})

The spectral theorem (See \cite[VI \S 5.3]{MR53:11389} ) asserts in our case that 
$$-L_\G=\sum_{i=0}^\infty\la_i P_i$$ where $P_i$ is the projection on the line spanned by $\psi_i$. We shall use $P_{\la_i}$ to denote the projection to the $\la_i$-eigenspace.

It is convenient to introduce a variable $s$ subject to the condition $\la=s(1-s)$. Hence the $s$-plane is a two-sheeted covering of the $\la$-plane and the  right half plane $\Re(s)>1/2$ cut along $1/2<s\leq 1$ corresponds to the $\la$-plane cut along the positive real axis.

The resolvent is the bounded operator $$R(s)=(L_{\G}+s(1-s))^{-1}:\L^2(\GmodH)\to \L^2(\GmodH)$$ defined for $s(1-s)\neq\la_i$. It satisfies 
\begin{equation}\label{resolventbound}\norm{R(s)}_\infty\leq \frac{1}{\min_{i}\abs{s(1-s)-\la_i}}\leq \frac{1}{\Im (s(1-s))}\leq\frac{1}{\abs{t}(2\s-1)} \end{equation}

From the spectral theorem we may conclude (See \cite[VI\S 5.2]{MR53:11389}) that
$$R(s)=\sum_{i=0}^\infty \frac{1}{s(1-s)-\la_i}P_i.$$ We note that for $s(1-s)$ close to $\la_i$ $$R(s)-\frac{1}{s(1-s)-\la_i}P_{\la_i}$$ is regular in $s$. Hence if $u(z,s)\in \L^2(\GmodH)$ is meromorphic in $s$ with a pole of order $k-1$ at $s_i$ with $s_i(1-s_i)=\la_i$ and leading term $u_{k-1}(z)$, then $R(s)u(z,s)$ has a pole of at most order $k$. The pole order is $k$ if and only if  
$$\lim_{s\to s_i}(s-s_i)^{k}R(s)u(z,s)=\frac{1}{1-2s_i}P_{\la_i}u_{k-1}(z)$$ is nonzero, and if this is the case then this is the leading term. If this is not the case the pole order is strictly less than $k$. 
We shall often use the above expression for $s_0=1$. Since $\psi_0=\vol(\GmodH)^{-1}$ this reduces to \begin{equation}\frac{-1}{\vol(\GmodH)}\int_{\GmodH}u_{k-1}(z)d\mu(z).\end{equation} 

\section{Hyperbolic Eisenstein series}
In this section we define hyperbolic Eisenstein series related to $\g_1$. The construction is a weight 0 real-analytic analogue of the holomorphic hyperbolic Eisenstein series of weight  $k\geq2$ considered in e.g. \cite{petersson,kudlamillson}. We shall only develop the theory of these hyperbolic Eisenstein series to the point needed to prove Theorem \ref{maintheoremforcomplex}. We shall have more to say about these series in \cite{hyperbolicrisager}.  We fix, in this section, a unitary character $\chi:\G\to S^1$ which is trivial on $\G_1$.

\begin{definition}The hyperbolic Eisenstein series related to $\g_1$ is defined by 
\begin{equation}
  E^{\g_1}(z,s)=\sum_{\g\in\GmodG}\overline{\chi(\g)}\left( \frac{\Im(\g z)}{\abs{\g z}}\right)^s
\end{equation} in its domain of absolute convergence. 
\end{definition}
It is easy to see that this is well defined in the domain of absolute convergence, and that it is $(\chi,\G)$ automorphic i.e. $$E^{\g_1}(\g z,s)=\chi(\g) E^{\g_1}(z,s)$$ whenever $\g\in\G$.

 \begin{lemma}\label{convergence} The series defining the hyperbolic Eisenstein series is absolutely convergent for $\Re(s)>1$. For $\Re(s)\geq l>1$ it is uniformly convergent and for  $(z,s)\in\H\times\{s\in \C\vert \Re(s)\geq l>1\}$ it is  bounded.
\end{lemma}
\begin{proof}The proof given is closely modeled after the proof of the convergence of the usual Eisenstein series given in \cite[Theorem 2.1.1]{kubota}. 
 We note that $$A=\{z\in \H\vert 1<\abs{z}<\mu \}$$ is a fundamental domain for $\G_1$. By the $\G_1$ invariance of $\Im(z)/\abs{z}$ we may assume that $\g z\in\overline{A}$ for all $\g\in \GmodG$. 

For any $\e>0$ we let $$k_\e(z,z')=\begin{cases}1 &\textrm{ if }d(z,z')\leq \e\\0&\textrm{ otherwise.}\end{cases}$$ Here $d(z,z')$ denotes the hyperbolic distance between $z$ and $z'$. As in \cite{kubota} we find that there exist $\Lambda_\e>0$ only dependent of $\e$ such that
$$\int_\H k_\e(z_0,z')\Im( z')^sd\mu(z')=\Lambda_\e \Im( z_0)^s.$$
If we choose $\e$ small enough we may assume that $$B(\g z,\e)\cap B(\g' z,\e)=\emptyset$$ when $\g\neq \g' \mod \G_1$. Here $B(z,\e)=\{z'\in H\vert d(z,z')<\e \}$. Hence we have, with $\Re(s)=\s$ 
\begin{eqnarray*}
  \sum_{\g\in\GmodG}\left( \frac{\Im(\g z)}{\abs{\g z}}\right)^\s&\leq& \frac{1}{\Lambda_\e}\sum_{\g\in\GmodG}\int_\H k_\e(\g z,z')\Im( z')^\s d\mu(z')\\
  &\leq &\frac{1}{\Lambda_\e}\sum_{\g\in\GmodG}\int_{B(\g z,\e)} \Im( z')^\s d\mu(z')
\\
&\leq &\frac{1}{\Lambda_\e}\int_{\overline{A}}\Im( z')^\s d\mu(z')\\
&\leq &\frac{1}{\Lambda_\e}\mu\int_0^{\mu}y^{\s-2}dy\\
&= &\frac{\mu^\s}{ ( \s-1) \Lambda_\e }
\end{eqnarray*}
From this inequality all the claims of the lemma easily follow.
\end{proof}
We note that the above proof also applies when the group is cofinite.
\begin{lemma}\label{differentialequation}The hyperbolic Eisenstein series satisfies
  \begin{equation}
    (\Delta+s(1-s)) E^{\g_1}(z,s)=-s^2 E^{\g_1}(z,s+2)
  \end{equation}
\end{lemma}
\begin{proof}
We note that since  $\Delta$ commutes with the $\SL(\R)$ action it suffices to show that 
  $$\left(\Delta +s(1-s)\right)\left(\frac{y}{\abs{z}}\right)^s =-s^2 \left(\frac{y}{\abs{z}}\right)^{s+2}$$
which is elementary. We omit the details.
\end{proof}
\begin{theorem}\label{continuation}
  The function $E^{\g_1}(z,s)$ has meromorphic continuation to the whole $s$ plane. At a regular point, $s_0$, $E^{\g_1}(z,s)$ is square integrable on $\GmodH$. The poles are located at $-2n+s_j$ where $s_j(1-s_j)$ is an eigenvalue of the automorphic Laplacian and $n\in\N$. If $\chi=1$ the pole at $s=1$ is simple with residue $$\frac{2\log \mu}{\vol{(\GmodH)}}.$$
\end{theorem}
\begin{proof} From Lemma \ref{convergence} we get that $E^{\g_1}(z,s)\in \L^2(\GmodH,d\mu(z))$ when $\Re(s)>1$ and we can therefore apply the resolvent to expression in Lemma \ref{differentialequation}. We get \begin{equation}\label{thisgivescontinuation}E^{\g_1}(z,s)=R(s)(-s^2E^{\g_1}(z,s+2)).\end{equation}  Since $E^{\g_1}(z,s)$ is holomorphic for $\Re(s)>1$ this gives meromorphic continuation to $\Re(s)>-1$. with poles possible poles at $s_j(1-s_j)$. Once this has been established, Eq. (\ref{thisgivescontinuation}) gives continuation to $\Re(s)>-3$. Repeating this process we obtain meromorphic continuation to the whole $s$-plane.    

The pole order at $s=1$ follows from the discussion in the end of section \ref{resolvent} and this also gives the residue
\begin{align*}
\frac{-1}{\vol(\GmodH)}&\int_{\GmodH}-1^2E^{\g_1}(z,3)d\mu(z)\\&=\frac{1}{\vol(\GmodH)}\int_{1\leq\abs{z}\leq \mu} \left(\frac{y}{\abs{z}}\right)^3 d\mu(z)
\\&=\frac{1}{\vol(\GmodH)}\int_{1}^\mu\int_{0}^\pi \left(\frac{r\sin \theta}{r}\right)^3 \frac{1}{(r\sin \theta)^2}rdrd\theta\\
\\&=\frac{2\log \mu}{\vol(\GmodH)}.
\end{align*}
\end{proof}

\section{Hyperbolic Eisenstein series twisted with modular symbols}
In this section we shall introduce hyperbolic Eisenstein series twisted with modular symbols. The analytic properties of these functions contains the information that eventually will enable us to conclude Theorem \ref{maintheoremforcomplex}.

Whenever $g$ is a holomorphic or harmonic 1-form  on $\GmodH$ we shall write 
$$\modsym{\g}{g}=-2\pi i\int_{z_0}^{\g z_0}g,$$ where $\g\in \G.$ We shall call this the modular symbol related to $g$. It is easy to see that this is independent of the path chosen and also that it is independent of $z_0$. We shall sometimes write $\modsym{\g}{f}$ instead of $\modsym{\g}{f(z)dz}$.

Let $\w_k$, $k=1\ldots n$, be holomorphic 1-forms on $\GmodH$. We have $\w_k(z)=f_k(z)dz$ where $f_k:\H\to\H$ is a modular form of weight 2. Let $\a_k=\Re(\w_k)$ or $\a_k=\Im(\w_k)$. We define 
\begin{equation}
\chi_{\ev}(\gamma )=\exp \left(
\sum_{k=1}^{n}\e_k\modsym{\g}{\a_k} \right).
\end{equation} If we assume that $\modsym{\g_1}{\a_k}=0$ for $k=1\ldots n$ we may construct an associated family of hyperbolic Eisenstein series by setting 
\begin{equation}\label{stpraxedis}
E^{\g_1}(z, s,\ev)=\sum_{\gamma\in\GmodG}                    
\chi_{\ev}(\gamma ) \left( \frac{\Im(\g z)}{\abs{\g z}}\right)^s
\end{equation}

 We will use the following convention. A function with a subscript variable will denote the partial derivative of the function in that variable. We have 
\begin{equation}\label{yougotmail}E^{\g_1}_{\e_1,\ldots,\e_n}(z,s,\ov)=\sum_{\gamma\in\GmodG}                  \prod_{k=1}^n\modsym{\g}{\a_k}  \left( \frac{\Im(\g z)}{\abs{\g z}}\right)^s,\end{equation}
when the sum is absolutely convergent. This is analogous to the function introduced by Goldfeld in \cite{goldfeld2, goldfeld1}. We notice that these functions may be seen as the  coefficients in a power series expansion in $\ev$ of $E^{\g_1}(z,s,\ev)$ around the point $\ev=\ov$.  Our aim in this chapter is to understand the analytic properties of this series in a neighborhood of the point $s=1$. It is these properties that will enable us later  to prove the distribution result stated in the introduction.

 We consider the space $\L^2(\GmodH , \bar\chi_{\ev})$ of square integrable func\-tions that transform as 
$$h(\gamma  z)=\bar\chi_{\ev}(\gamma )h(z), \quad \gamma\in\Gamma$$
under the action of the group. We introduce  unitary operators
$$U(\ev):\L^2(\Gamma\setminus\H )\to \L^2(\Gamma \setminus \H ,
\bar\chi_{\ev})$$
given by
$$(U(\ev)h)  (z):=U(z,\ev)h(z)=\exp\left(2\pi i \left(\sum_{k=1}^{n}\e_k \int_{z_0}^z\a_k \right)\right)h(z).$$

We set 
$$L(\ev)=U(\ev)^{-1}\Delta U(\ev)$$ and 
\begin{equation}\label{bananaboat}D^{\g_1}(z,s,\ev)=U(-\ev )E^{\g_1}(z,s,\ev).\end{equation}
Then using Lemma \ref{differentialequation} we see that  
\begin{equation}\label{differentiatingd}(L(\ev)+s(1-s))D^{\g_1}(z,s,\ev)=-s^2D^{\g_1}(z,s+2,\ev)\end{equation}

From (\ref{bananaboat}) we see that  $$D^{\g_1}(z,s,\ev)=\sum_{\g\in\GmodG}\exp\left(-2\pi i \left(\sum_{k=1}^{n}\e_k \int_{z_0}^{\g z} \a_k \right)\right)\left(\frac{\Im(\g z)}{\abs{\g z }}\right)^s.$$

By termwise differentiation  we find 
$$D^{\g_1}_{\e_1,\ldots,\e_n}(z,s,\ov)=\sum_{\g\in\GmodG}\prod_{k=1}^n\left(-2\pi i \int_{z_0}^{\g z} \a_k\right) \left(\frac{\Im(\g z)}{\abs{\g z}}\right)^s,$$ 
whenever the sum is absolutely convergent. (See Lemma \ref{squareintegrable} below).
From (\ref{bananaboat}) we also infer that
\begin{equation}\label{switchfunction1}D^{\g_1}_{\e_1,\ldots,\e_n}(z,s,\ov)=\!\!\!\!\sum_{\!\!\!\vec m\in\{0,1\}^n}\prod_{k=1}^n\left(-2\pi i\int_{z_0}^z \a_k\right)^{m_k}E^{\g_1}_{\e_1^{1-m_1},\ldots,\e_n^{1-m_n}}(z,s,\ov).\end{equation} 
\begin{equation}\label{switchfunction2}E^{\g_1}_{\e_1,\ldots,\e_n}(z,s,\ov)=\!\!\!\!\sum_{\!\!\!\vec m\in\{0,1\}^n}\prod_{k=1}^n\left(2\pi i\int_{z_0}^z\a_k\right)^{m_k}D^{\g_1}_{\e_1^{1-m_1},\ldots,\e_n^{1-m_n}}(z,s,\ov).\end{equation} 
These  relations between the $D^{\g_1}_{\e_1,\ldots,\e_n}(z,s,\ov)$ and the $E^{\g_1}_{\e_1,\ldots,\e_n}(z,s,\ov)$ functions  show that one family determines the other.

In order to ensure that the functions we shall be studying are well defined we need the following crude bound on the antiderivative of a  modular form of weight 2. We notice that the proof uses the same starting point as Hecke's bound on the Fourier coefficients of modular forms (See e.g. \cite[Chapter VII Theorem 4.5]{serre}) .
\begin{lemma}\label{crudebound} Let $z\in\{z\in\H\vert 1\leq \abs{z}\leq \mu\}$ and assume $\abs{z_0}>\mu$. Let $f$ be a modular form of weight 2. Then there exist $K\in\R_+$ such that 
\begin{equation*}\abs{\int_{z_0}^zf(w)dw}\leq K\frac{\abs{z}}{y}.\end{equation*}
\end{lemma}
\begin{proof} Since $f(\g z)=(cz+d)^2f(z)$ for $\g\in\G$ the function $|f(z)y|$ is $\G$-invariant. Since $\GmodH$ is compact it is bounded on $H$ i.e.
$$\abs{f(z)}\leq \frac{M}{y}$$ for some $M\in\R_+$. 

The integral in question is independent of the path chosen so we choose the direct line between $z_0$ and $z$. On the assumptions of the lemma we have $\abs{f(w)}\leq M/\Im(z)$ when $w$ is on the line between $z_0$ and $z$. Hence 
$$ \abs{\int_{z_0}^zf(w)dw}\leq M \frac{\abs{z-z_0}}{y}\leq M \frac{\abs{z}+\abs{z_0}}{y}\leq M(1+\abs{z_0}) \frac{\abs{z}}{y}.$$
\end{proof}

 The above lemma enables us to prove the following:
\begin{lemma}\label{squareintegrable} Let $n\geq 1$. For $\Re(s)$ sufficiently large we have 
  \begin{equation*}
    D^{\g_1}_{\e_1,\ldots,\e_n}(z,s,\ov)\in \L^2(\GmodH,d\mu).
  \end{equation*}
\end{lemma}
\begin{proof}
 Since $A=\{z\in\H\vert 1<\abs{z}<\mu \}$ is a fundamental domain for $\G_1$, we may choose representatives of $\GmodG$ such that $\g z$ is in $\overline{A}$. If we assume $\abs{z_0}>\mu$, we can now use Lemma \ref{crudebound} to conclude 
$$\abs{D^{\g_1}_{\e_1,\ldots,\e_n}(z,s,\ov)}\leq K \sum_{\g\in\GmodG}\left(\frac{\Im(\g z)}{\abs{\g z}}\right)^{s-n}.$$
The Lemma now follows from Lemma \ref{convergence}.  We may remove the assumption $\abs{z_0}>\mu$ by using (\ref{switchfunction1}) and (\ref{switchfunction2}). This is done by noticing that what we have proved already enables us to prove  that the sum representations of $E^{\g_1}_{\e_1^{1-m_1},\ldots,\e_n^{1-m_n}}(z,s,\ov)$ are convergent for $\Re(s)$ sufficiently large. Now these sums do not depend on $z_0$ and hence (\ref{switchfunction1}) shows that  $D^{\g_1}_{\e_1,\ldots,\e_n}(z,s,\ov)\in \L^2(\GmodH,d\mu)$ without the assumption on $\abs{z_0}$.
\end{proof}

We define  \begin{eqnarray*}\langle f_1 dz+f_2d\bar z, g_1dz+g_2d\bar z\rangle &=&2y^2(f_1\bar g_1+f_2\bar g_2)\\ 
\delta (pdx+qdy)&=&
-y^2(p_x+q_y).\end{eqnarray*}
\begin{lemma}\label{hetstraatje}
The conjugated operator $L(\ev) $ is given by
\begin{eqnarray}\label{viewofdelft}\nonumber
L(\ev)h&=&\Delta h+4\pi i \sum_{k=1}^{n}\epsilon_k\langle dh, \a_k\rangle-
 4\pi ^2 \left(\sum_{k,l=1}^{n}\e_k\e_l\modsym{\a_k}{\a_l}\right)h.
\end{eqnarray}
\end{lemma}
\begin{proof}The proof of  Lemma 2.2 in \cite{petridisrisager} may be used without changes. We notice that in the present case $\delta(\a_k)=0$. 
\end{proof}
Lemma \ref{hetstraatje} gives
\begin{equation}\label{Lderivedonce}
L_{\e_k}(\ov)h=4\pi i \langle dh, \a_k\rangle\end{equation}
\begin {equation}\label{Lderivedtwice}
L_{\e_k\e_l}(\ov)h=-8\pi ^2\langle \a_k, \a_l\rangle h .
\end{equation}
and all higher order derivatives vanish. Differentiating (\ref{differentiatingd})  
we get 
\begin{equation}\label{goingonce}
(\Delta+s(1-s))D^{\g_1}_{\e_k}(z, s,\ov)=-\left(L_{\e_k}(\ov)D^{\g_1}(z,s,\ov)\right)-s^2D^{\g_1}_{\e_k}(z,s+2,\ov)
\end{equation}
and
\begin{align}\label{goingtwice}
(\Delta+&s(1-s))D^{\g_1}_{\e_1,\ldots,\e_n}(z, s,\ov)=
-\left(\sum_{k=1}^{n}L_{\e_k}(\ov)D^{\g_1}_{\e_1,.,\widehat{\e_k},.,\e_n}(z,s,\ov)\right.\\\nonumber&\left.+\sum_{\substack{k,l=1\\ k<l}}^{n}L_{\e_k\e_l}(\ov)D^{\g_1}_{\e_1,.,\widehat{\e_k},.,\widehat{\e_l},.,\e_n}(z,s,\ov)\right)-s^2D^{\g_1}_{\e_1,\ldots,\e_n}(z,s+2,\ov).
\end{align}
 Here $\widehat{\e_k}$ means that we have excluded $\e_k$ from the list. When $\Re(s)$ is sufficiently large we can use Lemma \ref{squareintegrable} and invert (\ref{goingonce}) and (\ref{goingtwice}) by applying the resolvent of the automorphic  Laplace operator $R(s)=(\Delta_\Gamma +s(1-s))^{-1}$.  We get
\begin{equation}\label{denmarkiscold}
D^{\g_1}_{\epsilon_k}(z, s,\ov)=-R(s)\left(L_{\e_k}(\ov)D^{\g_1}(z,s,\ov)+s^2D^{\g_1}_{\e_k}(z,s+2,\ov)\right)
\end{equation}
and
\begin{align}\label{denmarkisrainfull}
D&^{\g_1}_{\e_1,\ldots,\e_n}(z, s,\ov)=-R(s)\left(\sum_{k=1}^{n}L_{\e_k}(\ov)D^{\g_1}_{\e_1,.,\widehat{\e_k},.,\e_n}(z,s,\ov)\right.\\\nonumber&\left.+\sum_{\substack{k,l=1\\ k<l}}^{n}L_{\e_k\e_l}(\ov)D^{\g_1}_{\e_1,.,\widehat{\e_k},.,\widehat{\e_l},.,\e_n}(z,s,\ov)+s^2D^{\g_1}_{\e_1,\ldots,\e_n}(z,s+2,\ov)\right).
\end{align}
 This will turn out to be identities of great importance for the proofs of many results in this and the following chapter. As a starting point we prove 
\begin{lemma}\label{applejack}The function $D^{\g_1}_{\e_1\ldots\e_n}(z, s,\ov)$ has meromorphic continuation to $s\in\C$. At a point of  regularity, $s_0$,  $D^{\g_1}_{\e_1\ldots\e_n}(z, s_0,\ov)$ is square integrable. In $\Re(s)>1$ the function is analytic.
\end{lemma}
\begin{proof} This is induction in $n$. For $n=0$ we quote Theorem \ref{continuation}, while the induction step follows from  (\ref{denmarkiscold}) and (\ref{denmarkisrainfull}).
\end{proof}
Using (\ref{switchfunction2}) and the above we  get the following theorem:
\begin{theorem}\label{twistedcontinuation}The function $E^{\g_1}_{\e_1\ldots\e_n}(z, s,\ov)$ has meromorphic continuation to $s\in\C$.  In $\Re(s)>1$ the function is analytic.
\end{theorem} This proves Theorem \ref{whattodo}.
\begin{corollary}\label{absconvergence}The sum (\ref{yougotmail}) is absolutely convergent for $\Re(s)>1$. 
\end{corollary}
\begin{proof}
 This follows from  a classical theorem due to Landau (see e.g \cite[Chapter VI, Proposition 2.7]{serre}).
\end{proof}
\begin{corollary}\label{nicebound}Let $f(z)dz$ be a holomorphic 1-form on $\GmodH$ such that  $\modsym{\g_1}{f}=0$
 For any fixed $z\in \H$, $\varepsilon>0$ we have 
\begin{equation*}
\int_{z_0}^{\gamma z_0}f(z)dz=o((\abs{az+b}\abs{cz+d})^\varepsilon)
\end{equation*} as $\abs{az+b}\abs{cz+d}\to \infty$. 
\end{corollary}
\begin{proof}
From Theorem \ref{twistedcontinuation} and (\ref{yougotmail}) one easily finds that for any $m\in\N$
$$\sum_{\g\in\GmodG}\left(\int_{z_0}^{\g z_0}f(z)dz\right)^m\left( \frac{\Im(\g z)}{\abs{\g z}}\right)^s$$ 
has meromorphic continuation to $\C$ and that it is is  analytic in $\Re(s)>1$. Using Landau's result again one may conclude that the above series is absolutely convergent for $\Re(s)>1$. Since the terms in an absolutely convergent series tends to zero we get that
 $$\left(\int_{z_0}^{\g z_0}f(z)dz\right)^m\left( \frac{\Im(\g z)}{\abs{\g z}}\right)^2=\left(\int_{z_0}^{\g z_0}f(z)dz\right)^m \frac{y^2}{\abs{cz+d}^2\abs{az+b}^2}$$ tends to zero as  $\abs{az+b}\abs{cz+d}\to \infty$. Hence 
$$\int_{z_0}^{\g z_0}f(z)dz=o(\abs{cz+d}^{2/m}\abs{az+b}^{2/m})$$ for any $m\in\N$. 
\end{proof}
We note that putting $z=i$ we obtain Theorem \ref{rathergoodbound}. 

We shall now show how we can obtain the  Laurent expansions of $D^{\g_1}_{\e_1,\ldots,\e_n}(z,s,\ov)$ from  (\ref{denmarkiscold}) and (\ref{denmarkisrainfull}).
We start by showing that $R(s)L_{\e_k}(\ov)D^{\g_1}_{\e_1,.,\widehat{\e_k},.,\e_n}(z,s,\ov)$ is regular. To this end we need the following lemma: 
\begin{lemma}\label{regular}
For all $j=1,\ldots,n$ and $\Re (s)>1$
  \begin{equation*}\int_{\GmodH}\langle dD^{\g_1}_{\e_1,\ldots,{\widehat{ \e_j}} ,\ldots\e_n}(z,s,\ov),\a_{j}\rangle d\mu(z)=0.\end{equation*}
 \end{lemma}
\begin{proof}
 We let $F$ be a fundamental polygon for the group $\G$. The domain $F$ is then an open connected set with a finite number of sides which are pairwise conjugated. We have
\begin{align}\nonumber\int_{F}&\langle d  D^{\g_1}_{\e_1,\ldots,\hat\e_j,\ldots\e_n}(z,s,\ov),\alpha_{j}\rangle d\mu(z)
=\\&\int_{F}\frac{\partial}{\partial z}\left( D^{\g_1}_{\e_1,\ldots,\hat\e_j,\ldots,\e_n}(z,s,\ov)\right)\frac{\overline{f_{j}}}{2}dxdy\label{twointegrals}\\\nonumber+&\int_{F}\frac{\partial}{\partial \overline{z}}\left(D^{\g_1}_{\e_1,\ldots,\hat\e_j,\ldots,\e_n}(z,s,\ov)\right)\frac{{f_{j}}}{2}dxdy.\end{align}

For any real-differentiable function $h:U\to\C$ where $U\subset \C$ and  any bounded domain $R\subset U$ with piecewise differentiable boundary Stokes theorem implies that \begin{equation*}2i\int_{R}\frac{\partial}{\partial \overline{z}}hdxdy=\int_{\partial R}hdz.\end{equation*} We apply this to the second integral in (\ref{twointegrals}).  Since $f_{j}$ is holomorphic, the integral equals
\begin{equation*}-\frac{i}{4}\int_{\partial F}D^{\g_1}_{\e_1,\ldots,\hat\e_j,\ldots,\e_n}(z,s,\ov){f_{j}}dz.\end{equation*} The boundary of the fundamental domain is the union of conjugated sides. These conjugated sides cancel in the integral and we get the result. 
\end{proof}

Using this we can now prove 
\begin{lemma} \label{regularfunction} The function \begin{equation}-R(s)L_{\e_j}(\ov)D_{\e_1\ldots,\widehat{\e_j},\ldots,\e_n}(z,s,\ov)\end{equation}  is regular at $s=1$.
\end{lemma}
\begin{proof} We shall write $B(z,s)=D_{\e_1,\ldots,\hat\e_j,\ldots\e_n}(z,s,\ov)$. 
From (\ref{denmarkisrainfull}) it is clear that $s=1$ is not an essential singularity. Assume that it is a pole of order $k>0.$ Hence
\begin{align}\label{contradictioncandidate}\lim_{s\to 1}(s-1)^k(-R(s)L_{\e_j}(\ov)B(z,s))\neq0.\end{align}
But
\begin{align*}\lim&_{s\to 1}(s-1)^k(-R(s)L_{\e_j}(\ov)B(z,s))\\
&=\lim_{s\to 1}(s-1)^{k-1}\left(\frac{1}{\vol(\GmodH)}\int_{\GmodH}L_{\e_j}(\ov)B(z',s)d\mu(z')\right)=0
\end{align*}
by  Section \ref{resolvent}, (\ref{Lderivedonce}) and Lemma \ref{regular}. This contradicts (\ref{contradictioncandidate}), which completes the proof.
\end{proof}
We let $\widetilde\Sigma_{2m}$ be the elements of the symmetric group on $2m$ 
letters $1,2,\ldots,2m$ for which $\sigma(2j-1)<\sigma(2j)$ for $j=1,\ldots,m$. We notice that this has $(2m)!/2^m$ elements, which is easily seen by induction. 
\begin{lemma}\label{Dleadingterms} 
If $n$ is even $D^{\g_1}_{\e_1,\ldots\e_n}(z,s,\ov)$ has a pole at $s=1$ of 
at most order $n/2+1$. The $(s-1)^{n/2+1}$ coefficient in the expansion of the function $D^{\g_1}_{\e_1,\ldots\e_n}(z,s,\ov)$ around $s=1$ is 
\begin{equation*} 
\frac{(-8\pi^2)^{n/2}2\log \mu}{\vol(\GmodH)^{n/2+1}}\sum_{\sigma\in \widetilde{\Sigma}_n}
\left(\prod_{r=1}^{n/2}\int_{\GmodH}\inprod{\alpha_{\sigma(2r-1)}}{\alpha_{\sigma(2r)}}d\mu(z)\right).
\end{equation*}
If $n$ is odd,  $D^{\g_1}_{\e_1,\ldots\e_n}(z,s,\ov)$ has a pole at $s=1$ of at most order $(n+1)/2$.
\end{lemma} 
\begin{proof} For $n=0$ we quote Lemma \ref{continuation}, and for $n=1$ (\ref{denmarkiscold}),  
Lemma \ref{regularfunction} and the discussion in the end of section \ref{resolvent} give the result. Assume that the result is true 
for all $n\leq n_0$. By (\ref{denmarkisrainfull}), (\ref{Lderivedtwice}), 
Lemma \ref{regularfunction} and the fact that $$(-R(s)(-8\pi^2\modsym{w_k}{w_l}D^{\g_1}_{\e_1,.,\hat\e_k,.,
\hat\e_l,.,\e_n}(z,s,\ov)))$$ can have pole order at most 1 more than $D^{\g_1}_{\e_1,.,\hat\e_k,.,
\hat\e_l,.,\e_n}(z,s,\ov))$ at $s=1$, we obtain the result about the pole orders. Note also that $$R(s)s^2D^{\g_1}_{\e_1,\ldots,\e_n}(z,s+2,\ov)$$ always contributes with at most a simple pole.
For even  $n$ we notice that by induction and using the discussion in the end of section \ref{resolvent}  we find that the $(s-1)^{n/2-1}$ coefficient is 
\begin{align*}\frac{-8\pi^2}{\vol{(\GmodH)}}&\frac{(-8\pi^2)^{(n-2)/2}\log \mu}{\vol(\GmodH)^{(n-2)/2+1}}\cdot\\\sum_{\substack{k,l=1\\ k<l}}^{n}\sum_{\sigma\in \widetilde{\Sigma}_{n-2}}&\left(\prod_{r=1}^{(n-2)/2}\hbox{}^\prime\int_{\GmodH}\!\!\!\!\!\inprod{\alpha_{\sigma(2r-1)}}{\alpha_{\sigma(2r)}}d\mu(z)\right)\int_{\GmodH}\!\!\!\!\!\inprod{\alpha_k}{\alpha_l}d\mu(z),\end{align*}
  where the prime in  the product means that we have excluded $\a_k,\a_l$ from the product and enumerated the remaining differentials accordingly. The result follows.
\end{proof}
Using this we can prove
\begin{theorem}\label{leadingterms} For all $n$  $E^{\g_1}_{\e_1,\ldots\e_n}(z,s,\ov)$ has a pole at $s=1$ of at most order $[n/2]+1$. If $n$ is even the  $(s-1)^{[n/2]+1}$ coefficient in the Laurent expansion of $E^{\g_1}_{\e_1,\ldots\e_n}(z,s,\ov)$ is \begin{equation*} \frac{(-8\pi^2)^{n/2}2\log \mu}{\vol(\GmodH)^{n/2+1}}\sum_{\sigma\in \widetilde{\Sigma}_n}\left(\prod_{r=1}^{n/2}\int_{\GmodH}\inprod{\alpha_{\sigma(2r-1)}}{\alpha_{\sigma(2r)}}d\mu(z)\right).\end{equation*}
\end{theorem}
\begin{proof}This follows from (\ref{switchfunction2}) and Lemma \ref{Dleadingterms}.\end{proof}

We notice that \begin{equation}\label{realreal}\modsym{\Re(f(z)dz)}{\Re(f(z)dz)}=\modsym{\Im(f(z)dz)}{\Im(f(z)dz)}=y^2\abs{f(z)}^2,\end{equation} while \begin{equation}\label{realimaginary}\modsym{\Re(f(z)dz)}{\Im(f(z)dz)}=0.\end{equation} Hence many of the involved integrals may be expressed in terms of the Petersson norm defined in the weight two case by  \begin{equation}\norm{f}=\left(\int_{\GmodH}y^2\abs{f(z)}^2d\mu(z)\right)^{1/2}.\end{equation} We shall write $E^{\Re^l,\Im^{n-l}}(z,s):=E_{\e_1,\ldots,\e_n}(z,s,\ov)$ where $\a_i=\Re(f(z)dz)$ for $i=1,\ldots,l$ and $\a_i=\Im(f(z)dz)$ for $i=l+1,\ldots,n$.

\section{Growth on vertical lines}
By Corollary \ref{absconvergence} we see that $E^{\g_1}_{\e_1,\ldots,\e_n}(z,s,\ov)=O_K(1)$ for $\Re(s)=\sigma>1$ and $z$ in a fixed compact set $K$. In this section we show that when we only require $\sigma>1/2$ then we still have at most polynomial growth on the line $\Re(s)=\sigma$.

We first prove : 
\begin{lemma}\label{reptilienzoo} The  hyperbolic Eisenstein series $E^{\g_1}(z, s)$ 
has polynomial growth in $s$ in $\Re (s)>1/2$. More precisely we have for any $\varepsilon>0$ and $1/2<\sigma\leq 1$
\begin{equation}
E^{\g_1}(z, \sigma+it)=O_K(\abs{t}^{6(1-\s)+\epsilon}).\end{equation}
\end{lemma}
\begin{proof}Fix $\s>1/2$ and let $s=\s+it$.  Using (\ref{thisgivescontinuation}), (\ref{resolventbound}) , lemmata \ref{convergence},\ref{differentialequation} and Theorem \ref{continuation} we get $\L^2(\GmodH)$-bounds
\begin{align*}\norm{E^{\g_1}(z,s)}_2&\leq \frac{\abs{s}^2}{\abs{t}(2\s-1)}\norm{E^{\g_1}(z,s+2)}_2\leq K_1\abs{t}\\
\norm{\Delta E^{\g_1}(z,s)}_2&=\norm{-s(1-s)E^{\g_1}(z,s)-s^2E^{\g_1}(z,s+2)}_2\\
&\leq\abs{s(1-s)}\norm{E^{\g_1}(z,s)}_2+\abs{s^2}\norm{E^{\g_1}(z,s+2)}_2\\
&\leq\left(\frac{\abs{s(1-s)}\abs{s^2}}{\abs{t}(2\s-1)}+\abs{s^2}\right)\norm{E^{\g_1}(z,s+2)}_2\leq K_2\abs{t}^3
\end{align*}
for $\abs{t}$ sufficiently large. To get a pointwise bound we use the Sobolev embedding theorem  (see \cite[6.22 Corollary (b)]{warner}). We denote by $\norm{\cdot}_{H^t}$ the Sobolev $t$-norm. So   in dimension $2$ the Sobolev embedding theorem implies that $$\norm{u}_{\infty}\leq c \norm{u}_{H^2},$$ where  for any  second order elliptic operator $P$ there exist a $c'$ such that $$\norm{u}_{H^2}\leq c'(\norm{u}_{2}+\norm{Pu}_{2}).$$ (See \cite[6.29]{warner}) We use $P=\Delta$ and get 
$$E^{\g_1}(z,s)=O(\abs{t}^3).$$ Applying Phragm\'en-Lindel\"of (see e.g. \cite[Appendix 5]{patterson})  in the strip $1/2+\d_1\leq \Re(s)\leq 1+\d_2$ gives the result.
\end{proof}
This proves the second part of Theorem \ref{imgoingtoaconcert}.
\begin{lemma}\label{reptilienzoo2}
The function $D^{\g_1}_{\e_1,\ldots,\e_n}(z, s,\ov)$ has polynomial growth in $t$ in $\Re (s)>1/2$. More precisely we have for any $\varepsilon>0$ and $1/2<\sigma\leq 1$ \begin{equation}
D^{\g_1}_{\e_1,\ldots,\e_n}(z, \sigma+it,\ov)
=O (\abs{t}^{(6(n+1))(1-\sigma)+\varepsilon}).\end{equation}
The involved constant depends on $\varepsilon$, $\sigma$ and $\a_1,\ldots,\a_n$ but not on $z$.
\end{lemma}
\begin{proof}
We note that it is enough to prove that for $\Re(s)>1/2$ we have $$D^{\g_1}_{\e_1,\ldots,\e_n}(z, \sigma+it,\ov)
=O (\abs{t}^{3(n+1)}).$$
Once this is established we can apply Phragmen-Lindel\"of in the strip $1/2+\d_1\leq \Re(s)\leq 1+\d_2$ to get the result.

We use induction in $n$. For $n=0$ we quote Lemma \ref{reptilienzoo} and note that we have exponent 3 without $\e$ in the proof of the lemma.
We now assume that for a fixed $\s>1/2$
\begin{eqnarray}\label{induction1}
D_{\e_1,\ldots,\e_m}(z, \sigma+it,\ov)
&=&O (\abs{t}^{3(m+1)})\\
\label{induction2}\norm{L_{\e_k}(\ov)D_{\e_1,.,\hat\e_k,.,\e_m}(z,s,\ov)}_2&=&O (\abs{t}^{3(m+1)})
\end{eqnarray}
whenever $m\leq n-1$. We want to give $\L^2(\GmodH)$-norm estimates on $D_{\e_1,\ldots,\e_m}(z, \sigma+it,\ov)$ and $\Delta D_{\e_1,\ldots,\e_m}(z, \sigma+it,\ov)$ so that we can  use the Sobolev embedding theorem as in the proof of Lemma \ref{reptilienzoo}. 

The identity which is going to boost the induction is (\ref{denmarkisrainfull}). 
We have, by the induction hypothesis and (\ref{Lderivedtwice}) the bound \begin{equation}\label{first}\norm{L_{\e_k\e_l}(\ov)D_{\e_1,.,\widehat{\e_k},.,\widehat{\e_l},.,\e_m}(z,s,\ov)}_2=O(\abs{t}^{3(n-1)})\end{equation}
We also have 
\begin{align}\label{teatime}\norm{L_{\e_k}(\ov)D^{\g_1}_{\e_1,.,\widehat{\e_k},.,\e_n}(z,s,\ov)}_2\leq C_1&\left(\norm{D^{\g_1}_{z,\e_1,.,\widehat{\e_k},.,\e_n}(z,s,\ov)}_2\right.\\
\nonumber&+\left.\norm{D^{\g_1}_{\overline{z},\e_1,.,\widehat{\e_k},.,\e_n}(z,s,\ov)}_2\right)
\end{align}
The first term can be estimated 
\begin{align*}
\norm{D^{\g_1}_{z,\e_1,.,\widehat{\e_k},.,\e_n}(z,s,\ov)}_2&\leq C_2\norm{D^{\g_1}_{\e_1,.,\widehat{\e_k},.,\e_n}(z,s,\ov)}_{H_1}\\
&\leq C_2\norm{D^{\g_1}_{\e_1,.,\widehat{\e_k},.,\e_n}(z,s,\ov)}_{H_2}\\
&\leq C_3\left(\norm{D^{\g_1}_{\e_1,.,\widehat{\e_k},.,\e_n}(z,s,\ov)}_{2}\right.\\&\qquad\left.+\norm{\Delta D^{\g_1}_{\e_1,.,\widehat{\e_k},.,\e_n}(z,s,\ov)}_{2}\right)
\end{align*}
The second term is $O(\abs{t}^{3n+2})$ which is seen from (\ref{goingtwice}) and the induction hypothesis. The first term is $O(\abs{t}^{3n})$ by the induction hypothesis.
The second term in (\ref{teatime}) can be handled in the same way so we conclude that 
\begin{equation}\label{second}\norm{L_{\e_k}(\ov)D^{\g_1}_{\e_1,.,\widehat{\e_k},.,\e_n}(z,s,\ov)}_2=O(\abs{t}^{3n+2}).\end{equation}
This certainly establishes (\ref{induction2}) when $m=n$.
We have \begin{equation}\label{third}\norm{-s^2D_{\e_1,\ldots,\e_n}(z,s+2,\ov)}_2=O(\abs{t}^2)\end{equation} by Lemma \ref{applejack}. By (\ref{denmarkisrainfull}), (\ref{resolventbound}) (\ref{first}), (\ref{second}) and (\ref{third}) we conclude 
\begin{equation}\label{tvison}\norm{D_{\e_1,\ldots,\e_m}(z, \sigma+it,\ov)}_2=O(\abs{t}^{3n+1}).\end{equation}
By (\ref{goingtwice}),  (\ref{first}), (\ref{tvison}), (\ref{second}) and (\ref{third}) we get 
\begin{equation}\label{tvisoff}\norm{\Delta D_{\e_1,\ldots,\e_m}(z, \sigma+it,\ov)}_2=O(\abs{t}^{3n+3}).\end{equation}
Using these two bounds we use the Sobolev embedding theorem as in the proof of Lemma \ref{reptilienzoo} we get 
$$\norm{D_{\e_1,\ldots,\e_m}(z, \sigma+it,\ov)}_\infty=O(\abs{t}^{3n+3})$$ which establish (\ref{induction2}) for $m=n$.
\end{proof}
Using the above lemma and (\ref{switchfunction2}) we conclude
\begin{theorem}\label{growthonvertical}
The function $E^{\g_1}_{\e_1,\ldots,\e_n}(z, s,\ov)$ has polynomial growth in $t$ for $\Re (s)>1/2$. More precisely we have for any $\varepsilon>0$, $1/2<\Re(s)\leq 1$ and $z$ in a compact set $K$ 
\begin{equation}
E^{\g_1}_{\e_1,\ldots,\e_n}(z, s,\ov)=O (\abs{t}^{(6(n+1))(1-\sigma)+\varepsilon}).\end{equation}
The involved constant depends on $\varepsilon$, $\sigma$, $K$ and $\a_1,\ldots,\a_n$.
\end{theorem} Hence we have proved at Theorem \ref{lasttheorem}.  
\section{Estimates of summatory functions}
In the two preceding sections we found the pole structure of the twisted hyperbolic Eisenstein series and we showed that as a function of $s$ this has at most polynomial growth on vertical lines. In this section we state and prove two technical propositions that enables us  to use these properties to get estimates on certain summatory functions.

We shall formulate the results in terms of general Dirichlet series. Fix $\{f_n\}\subset \R_+$ a non-decreasing series that tends to $\infty$ as $n\to\infty$. Let $a=\{a_n\}_{n=1}^\infty\subseteq \C$. For $s\in \C$ we consider \begin{equation}\label{generalsum}H^a(s)=\sum_{n=1}^\infty a_n f_n^{-s}\end{equation}
We assume that  
\renewcommand{\theenumi}{\roman{enumi}}
\begin{enumerate}
\item{\label{forste}The sum in (\ref{generalsum}) is absolutely convergent for $\Re(s)>1$}
\item{\label{anden}The function $H$ has meromorphic continuation to $\Re(s)>h-\varepsilon$ where $h<1$.}
\item{\label{tredje}The point $s=1$ is the only pole in $\Re(s)\geq h$.}
\item{\label{fjerde} The function  $H$ grows at most polynomially in $\Im(s)=t$ uniformly for  $h\leq \Re(s)\leq 2$, i.e. $H^a(s)=O(\abs{t}^b$)}
\end{enumerate}
We note that (\ref{forste}) implies that the series is uniformly convergent on compact subsets. By the Phragm\'en- Lindel\"of theorem we may replace (\ref{fjerde}) by the weaker assumption that for any fixed $h \leq \s\leq 2$, the function  $H^a(s)$ grows at most polynomially on vertical lines.
\begin{prop}\label{generaltheorem}Assume $\{a_n\}_{n=1}^\infty\subseteq \R_+$ and that $H^a(s)$ satisfies (\ref{forste})-(\ref{fjerde}). If $s=1$ is a simple pole then
\begin{equation}\sum_{f_n\leq T }a_n=\operatorname{Res}_{s=1}(H^a(s))T+O(T^{\frac{b+h}{b+1}+\varepsilon}).
\end{equation}
If $s=l$ is a pole of order $l>1$, and $d_l$ is the leading term in the Laurent expansion of $H^a(s)$ then 
\begin{equation}\sum_{f_n\leq T}a_n=\frac{d_l}{(l-1)!}T(\log T)^{l-1}+O(T(\log T)^{l-2}).
\end{equation}
\end{prop}
\begin{proof}
Let $\phi_U:\R\to\R$, $U\geq U_0$ ,  be a family of smooth decreasing functions with \begin{equation*}\phi_U(t)=\begin{cases}1 &\textrm{ if }t\leq1-1/U\\0 &\textrm{ if }t\leq1+1/U,  \end{cases}\end{equation*} and $\phi_U^{(j)}(t)=O(U^j)$ as $U\to\infty$. For $\Re(s)>0$ we let \begin{equation*}R_U(s)=\int_0^\infty\phi_U(t)t^{s-1}dt\end{equation*}be the Mellin transform of $\phi_U$. Then we have \begin{equation}\label{vud1}R_U(s)=\frac{1}{s}+O\left(\frac{1}{U}\right)\qquad\textrm{as }U\to\infty \end{equation} and for any $c>0$  \begin{equation}\label{vud2}R_U(s)=O\left(\frac{1}{\abs{s}}\left(\frac{U}{1+\abs{s}}\right)^c\right)\qquad\textrm{as }\abs{s}\to\infty.\end{equation} Both estimates are uniform for $\Re(s)$ bounded.   The first is a mean value estimate while the second is successive partial integration and a mean value estimate. 
The Mellin inversion formula now gives
\begin{align*}\sum_{n=1}^\infty a_n \phi_U\left(\frac{f_n}{T}\right)=&\sum_{n=1}^\infty a_n \frac{1}{2\pi i}\int_{\Re(s)=2}R_U(s)\left(\frac{f_n}{T}\right)^{-s}ds\\=&\frac{1}{2\pi i}\int_{\Re(s)=2}\!\!\!H^a(s)R_U(s)T^sds.
\end{align*} We note that by (\ref{vud2}) and (\ref{fjerde}) the integral is convergent. We now move the line of integration to the line $\Re(s)=h$ 
by integrating along a box of some height and then letting this height go to infinity. By (\ref{fjerde})  and   (\ref{vud2}) we find that the contribution from the horizontal sides goes to zero. Since we assume that $s=1$ is the only pole of the integrand with $\Re(s)\leq h$ then using Cauchy's residue theorem we obtain
\begin{align*}
\frac{1}{2\pi i}&\int_{\Re(s)=2}\!\!\! H^a(s)R_U(s)T^sds\\&=\res_{s=1}\left(H^a(s)R_U(s)T^s\right)+\frac{1}{2\pi i}\int_{\Re(s)=h}\!\!\! H^a(s)R_U(s)T^sds.
\end{align*}
If we choose $c=b+\varepsilon$ the last integral is convergent and $O(T^hU^{b+\varepsilon})$. 

Assume that $H^a(s)$ has a pole of order $l$ with $(s-1)^{-l}$ coefficient $d_{-l}$ then if $l>1$ we have
\begin{align*}
\res_{s=1}&\left(H^a(s) R_U(s)T^s\right)\\
&=\frac{1}{(l-1)!}\lim_{s\to 1}\frac{d^{l-1}}{ds^{l-1}}\left((s-1)^l\left(H^a(s)R_U(s)T^s\right)\right)\\
&=\frac{1}{(l-1)!}\!\!\!\!\sum_{n_1+n_2+n_3=l-1}\!\!\!\!\!\!\!\!\left.\frac{\partial^{n_1}(s-1)^l H^a(s)}{\partial s^{n_1}}\right|_{s=1}\!\!\!\!\left.\frac{\partial^{n_2}R_U(s) }{\partial s^{n_2}}\right|_{s=1}\!\!\!\!\left.\frac{\partial^{n_3} T^s}{\partial s^{n_3}}\right|_{s=1}\\
\intertext{The first factor in the sum is independent of $U$ and $T$, while the second is independent of $T$ and bounded in $U$. The third factor has leading term  $T(\log T)^{n_3}$ and a reminder $O(\log T^{n_3-1})$. Hence the leading term is the one corresponding to $n_1=n_2=0$, $n_3=l-1$ and we get, using (\ref{vud1}),}
&=\frac{d_{-l}}{(l-1)!y}T(\log T)^{l-1}+O(T(\log T)^{l-2}+T\log T^{l-1}/U). 
\end{align*}
This gives 
\begin{align*}\sum_{f_n\leq T}\omega_\gamma\phi_U\left(\frac{f_n}{T}\right)=& \frac{d_{-l}}{(l-1)!y}T(\log T)^{l-1}\\&+ O(T(\log T)^{l-2}+T\log T^{l-1}/U+T^hU^{a+\varepsilon}).\end{align*}
If $l=1$ then by (\ref{vud1})
\begin{equation*}
\res_{s=1}\left(H^a(s)R_U(s)T^s\right)=\frac{a_{-1}}{y}T+O(T/U),
\end{equation*}
and we get
\begin{equation*}
 \sum_{f_n\leq T}a_n \phi_U\left(\frac{f_n}{T}\right)= d_{-1}T+O(T/U+T^hU^{b+\varepsilon}).
\end{equation*}
 If $H^a(s)$ has a nonsimple pole we choose $U=\log T$ and we get
\begin{equation}\label{exp1}\sum_{f_n\leq T}a_n\phi_U\left(\frac{f_n}{T})\right)=\frac{d_{-l}}{(l-1)!}T(\log T)^{l-1}+O(T(\log T)^{l-2}).\end{equation}
In the simple pole case we choose $U=T^{(1-h)/(b+1+\varepsilon)}$ in order to balance the error terms and we get
\begin{equation}\label{exp2}\sum_{f_n\leq T} a_n\phi_U\left(\frac{f_n}{T})\right)=d_{-1}T+O(T^{\frac{b+h+\varepsilon}{b+1+\varepsilon}}).
\end{equation}  
At this point we note that if $a_n$ is non-negative for all $n$, then by further requiring $\phi_U(t)=0$ if $t\geq 1$ and  $\tilde\phi_U(t)=1$ for $t\leq 1$, we have 
\begin{equation*}\sum_{f_n\leq T} a_n \phi_U\left(\frac{f_n}{T}\right)\leq \sum_{f_n\leq T}a_n \leq \sum_{f_n\leq T} a_n \tilde\phi_U\left(\frac{f_n}{T}\right)
\end{equation*}
from which it easily follows that the middle sum has an asymptotic expansion.  In the case $s=1$ a simple pole we choose $U=T^\frac{1-h}{2-h+\varepsilon}$ to balance the error terms\end{proof}
Since $a_n\not\in\R$ for many of the applications we have in mind we shall also deal with this situation We let $H(s)$ be the sum corresponding to $a_n=1$ for all $n$. 
\begin{prop}\label{anothergeneraltheorem}Assume that $H(s)$ satisfies the conditions of Proposition \ref{generaltheorem} with parameters $h'$, $b'$. Assume that for any $\varepsilon>0$ we have $a_n=O((f_n)^\varepsilon)$ as $n\to\infty$ and that $H^a(s)$ satisfies (\ref{forste})-(\ref{fjerde}). Assume further that$$\ord_{s=1}(H^a(s))\geq\ord_{s=1}(H(s))$$ If $s=1$ is a simple pole then
\begin{equation}\sum_{f_n\leq T }a_n=\operatorname{Res}_{s=1}(H^a(s))T+O(T^{\max \left(\frac{b+h}{b+1},\frac{b'+h'}{b'+1}\right)+\varepsilon}).
\end{equation}
If $s=l$ is a pole of order $l>1$, and $c_l$ is the leading term in the Laurent expansion of $H^a(s)$ then 
\begin{equation}\sum_{f_n\leq T}a_n=d_lT(\log T)^{l-1}+O(T(\log T)^{l-2}).
\end{equation}
\end{prop}
\begin{proof}
We may re-use most of the proof of the last proposition. To get a result without $\phi_U$ from (\ref{exp1}) and (\ref{exp2}) we notice that if we choose $\phi_U$ such that $\phi_U(t)=1$ for $t\leq 1$ then 
\begin{equation*}
\sum_{f_n\leq T}a_n\phi_U\left((f_nT)^{-1}\right)=\sum_{f_n\leq T}a_n+\sum_{T<f_n\leq T(1+1/U)}a_n\phi_U\left((f_nT)^{-1}\right).
\end{equation*}
From  $a_n=O((f_n)^\varepsilon)$  we see that we may evaluate the last sum in the following way. For any $\varepsilon>0$ this is less than a constant times
\begin{equation*}(T(1+1/U))^\varepsilon\!\!\!\!\!\! \sum_{T<f_n\leq T(1+1/U)}\!\!\!\!\!\!1 \leq (2T)^\varepsilon \!\!\!\!\!\!\sum_{T<f_n\leq T(1+1/U)}\!\!\!\!\!\!1.
\end{equation*} 
By Proposition \ref{generaltheorem} the sum is $O(T/U)+O(T^{\frac{b'+h}{b'+1}+\varepsilon})$ if $H(s)$ has a simple pole and $O(T/U(\log T)^{l'-1})+O(T(\log T)^{l'-2})$ if $\ord_{s=1}(H(s))=l'>1$. Using this with the above choices of $U$ we get the result.
\end{proof} We note that under the assumptions of the above proposition, with  the exceptions  that $H^a(s)$should  be regular at $s=1$ and that $H(z,s)$ should have  a simple pole at $s=1$, we may conclude that $$\sum_{f_n\leq T }a_n=O(T^{\max \left(\frac{b+h}{b+1},\frac{b'+h'}{b'+1}\right)+\varepsilon}).$$ The proof of this is identical to the proof of the above with $d_{-1}=0$.

We now observe that by ordering the elements of $\GmodG$ by defining $\g\leq\tilde\g$ if and only if $$\frac{\Im(\tilde\g z)}{\abs{\tilde\g z}}\leq \frac{\Im(\g z)}{\abs{\g z}}$$ we get 
$$\GmodG=\left\{\g_2,\g_3,\ldots\g_n,\ldots\right\}$$ with $\g_n\leq\g_{n+1}$. If we define $$f_n=\left(\frac{\Im(\g_n z)}{\abs{\g_n z}}\right)^{-1}$$ then $f_n$ is a non-decreasing series tending to $\infty$ and we have $H(s)=E^{\g_1}(z,s)$. Observing that $$\frac{\Im(\g z)}{\abs{\g z}}=\frac{y}{\abs{az+b}\abs{cz+d}}$$ Theorem \ref{continuation}, Lemma \ref{reptilienzoo} and Proposition \ref{generaltheorem} enables us to conclude
\begin{theorem}\label{counting}
\begin{equation*}\sum_{\substack{\gamma\in\GmodG\\\abs{az+b}\abs{cz+d}\leq T}}1=\frac{2\log\mu}{y\vol(\GmodH)}T+O(T^{1-\delta})\end{equation*} for some $\delta>0$.  \end{theorem}
Setting $z=i$ we obtain Theorem \ref{homotopyclasses}.
The estimation of the remainder term depends on the growth estimates of $E^{\g_1}(z,s)$ available to us and the existence of small eigenvalues of the Laplacian. Assuming no small eigenvalues Lemma \ref{reptilienzoo} enables us to conclude $1-\delta=7/8+\varepsilon$.

We now fix a holomorphic 1-form $f(z)dz$ Let $\alpha=\Re(f(z)dz)$ and $\beta=\Im(f(z)dz)$.
\begin{theorem}\label{sum}We have 
\begin{align}
\label{evensum}\sum_{\substack{\gamma\in\GmodG\\\abs{az+b}\abs{cz+d}\leq T}}&\!\!\!\modsym{\gamma}{\alpha}^{2m}\modsym{\gamma}{\beta}^{2n}\\\nonumber=&\frac{(-8\pi^2)^{m+n}\norm{f}^{2m+2n}2\log \mu }{y\vol(\GmodH)^{m+n+1}}&\frac{(2m)!}{m!2^{m}}\frac{(2n)!}{n!2^{n}}T\log^{m+n} T\\\nonumber&&+O(T\log^{m+n-1} T),\end{align} 
and if $m$ or $n$ is odd then
\begin{align}\label{oddsum}\sum_{\substack{\gamma\in\GmodG\\\abs{az+b}\abs{cz+d}\leq T}}\modsym{\gamma}{\alpha}^{m}\modsym{\gamma}{\beta}^{n}&=O(T\log^{k} T) & 
\end{align} for some $k\in \N$ strictly less than $(m+n)/2$.
\end{theorem}
\begin{proof}
This follows from Proposition \ref{anothergeneraltheorem} with $a_k=\modsym{\g_k}{\a}^m\modsym{\g_k}{\beta}^n$. We have $a_n=O(f_n^
\varepsilon)$ by Corollary \ref{nicebound}. The assumptions of Proposition \ref{anothergeneraltheorem} are satisfied by Theorem \ref{continuation}, Lemma \ref{reptilienzoo} Theorem \ref{leadingterms} and Theorem \ref{growthonvertical}.  We also notice that \begin{equation}\frac{(2m)!(2n)!}{2^{m+n}(m+n)!}\binom{m+n}{n}=\frac{(2m)!}{m!2^{m}}\frac{(2n)!}{n!2^{n}}.\end{equation}
\end{proof}
\section{The distribution of  modular symbols}
We now show how to obtain a distribution result for the modular symbols from the asymptotic expansions of Corollary \ref{sum}. We renormalize the modular symbols in the following way. Let \begin{eqnarray*}\label{renormalize}\widetilde{\modsym{\gamma}{f}}&=&\sqrt{\frac{\vol{(\GmodH)}}{{8\pi^2\norm{f}^2}}}\modsym{\gamma}{f}\\
\widetilde{\modsym{\gamma}{\alpha}}&=&\sqrt{\frac{\vol{(\GmodH)}}{{8\pi^2\norm{f}^2}}}\modsym{\gamma}{\alpha}\\
\widetilde{\modsym{\gamma}{\beta}}&=&\sqrt{\frac{\vol{(\GmodH)}}{{8\pi^2\norm{f}^2}}}\modsym{\gamma}{\beta}\\
 \end{eqnarray*} where  $\alpha=\Re(f(z)dz)$, $\beta=\Im(f(z)dz)$. Let furthermore \begin{equation}(\GmodG)^T:=\left\{\gamma\in\GmodG |\quad \abs{az+b}\abs{cz+d}\leq T\right\}.\end{equation} By Theorem \ref{counting} we have \begin{equation}\label{countingexpr}\#(\GmodG)^T=\frac{2\log \mu T}{\vol({\GmodH})y}+O(T^{1-\delta}),\end{equation} for some $\delta>0$. 
Now let $X_T$ be the random variable with probability measure
\begin{equation}P(X_T\in R)=\frac{ \#\left\{\gamma\in (\GmodG)^T \left| \frac{\widetilde{\langle \gamma,f\rangle}}{\sqrt{\log \abs{az+b}\abs{cz+d}}}\in R\right. \right\} }{\#(\GmodG)^T}.
\end{equation} for $R\subset\C$ (By convention we set  $\widetilde{<\gamma,\alpha>}/{\sqrt{\log \abs{az+b}\abs{cz+d}}}=0$ if $\abs{az+b}\abs{cz+d}\leq 1$. Note that there are only finitely many such elements.) We consider the moments of $X_T$
\begin{equation}M_{n,m}(X_T)=\sum_{\gamma\in (\GmodG)^T}\frac{\left[\Re\left(\frac{\widetilde{\modsym{\gamma}{f}}}{\sqrt{\log \abs{az+b}\abs{cz+d}}}\right)\right]^n\left[\Im\left(\frac{\widetilde{\modsym{\gamma}{f}}}{\sqrt{\log \abs{az+b}\abs{cz+d}}}\right)\right]^m}{\#(\GmodG)^T},\end{equation}
and note that \begin{eqnarray*}
\Re(\widetilde{\modsym{\gamma}{f}})&=&i\widetilde{\modsym{\gamma}{\beta}}\\
\Im(\widetilde{\modsym{\gamma}{f}})&=&-i\widetilde{\modsym{\gamma}{\alpha}}.\\
\end{eqnarray*}  
By partial summation we have   
\begin{align*}M_{n,m}(X_T)=&\frac{i^{n+m}(-1)^m}{\#(\GmodG)^T}\left(\sum_{\gamma\in(\GmodG)^T}\widetilde{\modsym{\gamma}{\beta}}^n\widetilde{\modsym{\gamma}{\alpha}}^m\frac{1}{\log T^{(m+n)/2}}\right.\\
&\!\!\left.+\frac{m+n}{2}\!\!\int_{\delta_0}^T\!\!\sum_{\gamma\in(\GmodG)^t}\widetilde{\modsym{\gamma}{\beta}}^n\widetilde{\modsym{\gamma}{\alpha}}^m\frac{1}{t(\log t)^{(m+n)/2+1}}dt\right),\end{align*} where $\delta_0=\min\{\abs{az+b}\abs{cz+d}:\abs{az+b}\abs{cz+d}>1\}$. If we now apply Corollary \ref{sum} and (\ref{countingexpr}) we find that as $T\to\infty$
\begin{equation}M_{n,m}(X_T)\to \begin{cases} \frac{n!}{(n/2)!2^{n/2}}\frac{m!}{(m/2)!2^{m/2}},  &\textrm{ if }m\textrm{ and }n \textrm{ are even},\\0,  &\textrm{ otherwise. }\end{cases}\end{equation}
We notice that the right-hand side is the moments of the bivariate Gaussian  distribution with correlation coefficient zero. Hence by a result due to Fr\'echet and Shohat (see \cite[11.4.C]{loeve}) we conclude the following:
\begin{theorem}\label{main}Asymptotically  $\frac{\widetilde{\langle\gamma,f\rangle}}{\sqrt{\log \abs{az+b}\abs{cz+d}}}$ has bivariate Gaussian distribution with correlation coefficient zero. More precisely we have
\begin{equation}\frac{ \#\left\{\gamma\in (\GmodG)^T \left| \frac{\widetilde{\langle\gamma,f\rangle }}{\sqrt{\log \abs{az+b}\abs{cz+d}}}\in R\right. \right\} }{\#(\GmodG)^T}\to \frac{1}{{2\pi}}\int_R\!\! \exp\left(-\frac{x^2+y^2}{2}\right)dxdy\end{equation} as $T\to\infty$.
\end{theorem}
As an easy corollary we get the following result about the distribution of harmonic differentials
\begin{corollary}\label{mainreal}Asymptotically  $\frac{\Re(\widetilde{\langle \gamma,f\rangle })}{\sqrt{\log \abs{az+b}\abs{cz+d}}}$ has Gaussian distribution. More precisely we have
\begin{equation}\frac{ \#\left\{\gamma\in (\GmodG)^T \left| \frac{\Re(\widetilde{\langle\gamma,f\rangle )}}{\sqrt{\log \abs{az+b}\abs{cz+d}}}\in [a,b]\right. \right\} }{\#(\GmodG)^T}\to \frac{1}{\sqrt{2\pi}}\int_a^b\!\! \exp\left(-\frac{x^2}{2}\right)dx\end{equation} as $T\to\infty$.
\end{corollary}
The same holds for $ \Im(\widetilde{\modsym{\gamma}{f}})$. We note that by putting $z=i$ in Corollary \ref{mainreal} and Theorem \ref{main} we obtain  Theorem \ref{maintheoremforcomplex} and Theorem \ref{maintheoremforreal}.

\bibliographystyle{plain}

\begin{thebibliography}{10}

\bibitem{buser}
Peter Buser.
\newblock {\em Geometry and spectra of compact {R}iemann surfaces}, volume 106
  of {\em Progress in Mathematics}.
\newblock Birkh\"auser Boston Inc., Boston, MA, 1992.

\bibitem{goldfeld2}
Dorian Goldfeld.
\newblock The distribution of modular symbols.
\newblock In {\em Number theory in progress, Vol. 2 (Zakopane-Ko\'scielisko,
  1997)}, pages 849--865. de Gruyter, Berlin, 1999.

\bibitem{goldfeld1}
Dorian Goldfeld.
\newblock Zeta functions formed with modular symbols.
\newblock In {\em Automorphic forms, automorphic representations, and
  arithmetic (Fort Worth, TX, 1996)}, volume~66 of {\em Proc. Sympos. Pure
  Math.}, pages 111--121. Amer. Math. Soc., Providence, RI, 1999.

\bibitem{MR53:11389}
Tosio Kato.
\newblock {\em Perturbation theory for linear operators}.
\newblock Springer-Verlag, Berlin, second edition, 1976.
\newblock Grundlehren der Mathematischen Wissenschaften, Band 132.

\bibitem{kubota}
Tomio Kubota.
\newblock {\em Elementary theory of {E}isenstein series}.
\newblock Kodansha Ltd., Tokyo, 1973.

\bibitem{kudlamillson}
Stephen~S. Kudla and John~J. Millson.
\newblock Harmonic differentials and closed geodesics on a {R}iemann surface.
\newblock {\em Invent. Math.}, 54(3):193--211, 1979.

\bibitem{loeve}
Michel Lo{\`e}ve.
\newblock {\em Probability theory. {I}}.
\newblock Springer-Verlag, New York, fourth edition, 1977.
\newblock Graduate Texts in Mathematics, Vol. 45.

\bibitem{patterson}
S.~J. Patterson.
\newblock {\em An introduction to the theory of the {R}iemann zeta-function},
  volume~14 of {\em Cambridge Studies in Advanced Mathematics}.
\newblock Cambridge University Press, Cambridge, 1988.

\bibitem{petersson}
Hans Petersson.
\newblock Ein {S}ummationsverfahren f\"ur die {P}oincar\'eschen {R}eihen von
  der {D}imension --2 zu den hyperbolischen {F}ixpunktepaaren.
\newblock {\em Math. Z.}, 49:441--496, 1944.

\bibitem{petridisrisager}
Yiannis Petridis and Morten Skarsholm~Risager.
\newblock Modular symbols are normally distributed.
\newblock Submitted, 2003.

\bibitem{serre}
J.-P. Serre.
\newblock {\em A course in arithmetic}.
\newblock Springer-Verlag, New York, 1973.
\newblock Translated from the French, Graduate Texts in Mathematics, No. 7.

\bibitem{MR47:3318}
Goro Shimura.
\newblock {\em Introduction to the arithmetic theory of automorphic functions}.
\newblock Publications of the Mathematical Society of Japan, No. 11. Iwanami
  Shoten, Publishers, Tokyo, 1971.
\newblock Kan\^o Memorial Lectures, No. 1.

\bibitem{hyperbolicrisager}
Morten Skarsholm~Risager.
\newblock Hyperbolic {E}isenstein series.
\newblock In preparation.

\bibitem{warner}
Frank~W. Warner.
\newblock {\em Foundations of differentiable manifolds and {L}ie groups},
  volume~94 of {\em Graduate Texts in Mathematics}.
\newblock Springer-Verlag, New York, 1983.
\newblock Corrected reprint of the 1971 edition.

\end{thebibliography}

\end{document}